# LLM-Enhanced, Data-Driven Personalized and Equitable Clinician Scheduling: A Predict-then-Optimize Approach

Anjali Jha[1], Wanqing Chen[1], Maxim Eckmann[2], Ian Stockwell[1], Jianwu Wang[1], Kai Sun[1,2]
[1]*University of Maryland, Baltimore County*, Baltimore, MD, USA
[2]*University of Texas Health Science Center at San Antonio*, San Antonio, TX, USA
Email: [1]{qd98375, wchen7, istock1, jianwu, ksun1}@umbc.edu, [2]{eckmann, sunk1}@uthscsa.edu

*Abstract*—Clinician scheduling remains a persistent challenge due to limited clinical resources and fluctuating demands. This complexity is especially acute in large academic anesthesiology departments as physicians balance responsibilities across multiple clinical sites with conflicting priorities. Further, scheduling must account for individuals' clinical and lifestyle preferences to ensure job satisfaction and well-being. Traditional approaches, often based on statistical or rule-based optimization models, rely on structured data and explicit domain knowledge. However, these methods often overlook unstructured information, e.g., free-text notes from routinely administered clinician well-being surveys and scheduling platforms. These notes may reveal implicit constraints (post-call recovery needs, early departures for family obligations) and underutilized clinical resources (a physician's willingness to work additional hours to address emergent needs). Neglecting such information can lead to misaligned schedules, increased burnout, overlooked staffing flexibility, and suboptimal utilization of available resources. To address this gap, we propose a predict-then-optimize framework that integrates classification-based clinician availability predictions with a mixed-integer programming schedule optimization model. Large language models (LLMs) are employed to extract actionable preferences and implicit constraints from unstructured schedule notes, enhancing the reliability of availability predictions. These predictions then inform the schedule optimization considering four objectives: (i) ensuring clinical full-time equivalent compliance, (ii) reducing workload imbalances by enforcing equitable proportions of shift types, (iii) maximizing clinician availability for assigned shifts, and (iv) schedule consistency. By combining the interpretive power of LLMs with the rigor of mathematical optimization, our framework provides a robust, data-driven solution that enhances operational efficiency while supporting equity and clinician well-being.

*Keywords—large language models, predict-then-optimize, healthcare operations, mixed-integer programming, data-driven optimization, multi-objective optimization*

## I. Introduction

This paper addresses the problem of scheduling clinicians in an outpatient pain clinic, where long planning horizons, physician redeployments, and unanticipated absences make it difficult to generate fair and reliable schedules. These difficulties exemplify a broader challenge in clinician scheduling, one of the most complex problems in healthcare operations, where limited resources and fluctuating demand must be balanced against contractual obligations and individual preferences. In large academic anesthesiology departments, these pressures are especially acute, as clinicians juggle responsibilities across diverse clinical sites with competing priorities. Effective scheduling must not only ensure adequate coverage but also account for clinicians' preferences to support job satisfaction, retention, and overall well-being [1]. Failure to achieve this balance can result in misaligned schedules, reduced workforce morale, and higher risks of burnout.

Healthcare systems globally are facing mounting pressure due to rising demand, increasing costs, and changing reimbursement policies [2]. These pressures have intensified longstanding anesthesiologist and nursing shortages [3-5], making efficient scheduling tools essential to optimize the deployment of this scarce workforce [6]. Within this context, outpatient pain clinics affiliated with academic anesthesiology departments face unique scheduling challenges. Because pain management care is non-emergent [7], staffing for pain clinics is often deprioritized when anesthesiologists are needed in higher-acuity clinical sites, e.g., operating rooms (ORs) or intensive care units (ICUs). Additionally, the work schedules need to balance non-clinical duties, e.g., teaching, administration, and research [8]. Persistent shortages and rising demand exacerbate this imbalance, leading to frequent redeployments, appointment cancellations, and reduced patient access [9].

Further complexity arises from mismatched scheduling horizons, i.e., outpatient pain clinics typically plan schedules three months in advance, while OR and ICU schedules are set monthly. Even after schedules are published, clinicians may be redeployed to cover urgent needs elsewhere. Additional uncertainty arises from unanticipated absences, e.g., sick leave or family obligations [10]. Mitigation strategies, e.g., scheduling daily backup clinicians, often conflict with administrative duties or exceed contractual clinical full-time equivalent (cFTE) allocations, disrupting workload equity and further reducing satisfaction.

Current practice in many clinics reflects these limitations. Clinician preferences entered into scheduling platforms, e.g., Qgenda, are often reviewed manually, inconsistently applied, or ignored altogether. Even when considered, they are handled case by case rather than systematically integrated into scheduling decisions. This inconsistency undermines transparency and equity, particularly in assigning additional compensated or benefit-enhanced shifts, e.g., additional services coverage.

To better utilize limited resources in the face of fluctuating demand and clinical resource uncertainty, data-driven scheduling approaches have become essential [11, 12]. These approaches typically rely on structured data and explicit domain rules, often formulated as stochastic optimization or mixed-integer programming (MIP) models [10, 13]. While such models provide mathematical rigor and can enforce fairness constraints, e.g., cFTE allocations, equitable shift distributions, they frequently overlook unstructured information embedded in free-text clinician notes, wellness surveys, and scheduling platforms. These sources often contain valuable implicit constraints, e.g., post-call recovery needs, early departure requests, or willingness to cover additional hours. Neglecting such information reduces scheduling flexibility, leads to underutilized resources, and increases inequities in workload distribution.

Recent advances in machine learning (ML) have also been applied to clinical resource planning, primarily for forecasting clinical demand, e.g., patient length of stay, surgical volumes, emergency department arrivals, and intensive care unit occupancy. While these models provide valuable insights into expected workload, they are rarely integrated into downstream scheduling or staffing decisions. However, predictions remain descriptive rather than prescriptive, limiting their operational impact. Large language models (LLMs), despite their growing role in medical text analysis and documentation [18, 19, 20], remain underutilized in clinician scheduling. In particular, they are not yet widely employed to extract day-level availability constraints from scheduling notes for direct integration into optimization models.

The predict-then-optimize (PTO) paradigm (also called the decoupled learning approach [21]) directly links predictions to decisions. It first predicts the unknown parameters of an optimization model with an ML-based predictor and then embeds those predictions in an optimization framework to produce actionable schedules. This paradigm provides a natural pathway for addressing outpatient clinician scheduling challenges.

Building on prior work in equitable anesthesiologist scheduling under demand uncertainty [22], we propose a two-step PTO clinician scheduling (PTO-CS) framework. **Step 1** is a classification-based clinician availability prediction model enhanced by LLMs to better extract latent and personalized constraints and preferences from both structured and unstructured historical scheduling data (_Contribution 1_). **Step 2** embeds these predictions into a multi-objective MIP that (i) ensures compliance with cFTE allocations, (ii) reduces workload imbalances by enforcing equitable proportions of shift types, (iii) maximizes clinician availability for assigned shifts, and (iv) maintains schedule consistency (_Contribution 2_). By combining LLM-enhanced availability predictions with the prescriptive rigor of optimization, our framework provides a robust and scalable solution to improve efficiency, equity, and clinician satisfaction.

To guide the reader, the remainder of this paper is organized as follows. Section II introduces the proposed PTO-CS framework, detailing the data sources, predictive modeling pipeline, and optimization formulation. Section III describes the experimental design and retrospective evaluation, including setup, metrics, and solver environment. Section IV presents the results, analyzing distributional deviations, match accuracy, and fairness outcomes across six months of clinic data. Section V concludes with key insights and discusses limitations and future directions.

## II. THE PREDICT-THEN-OPTIMIZE PARADIGM

The PTO paradigm links predictive modeling with prescriptive optimization in a two-step pipeline (see Fig. 1). In **Step 1**, an ML model is trained to predict uncertain parameters of an optimization problem. In **Step 2**, these predictions are embedded into a mathematical programming model that generates actionable decisions. This structure allows data-driven learning to inform operational optimization.

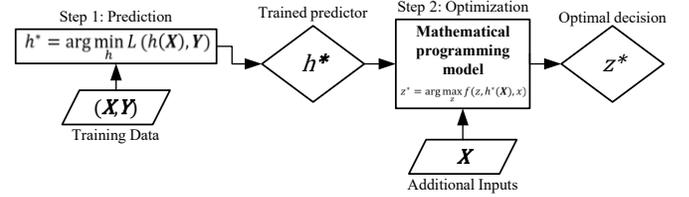

Fig. 1. A generic PTO paradiam

Formally, let $X$ denote the feature set for each clinician–day pair. In our clinician scheduling application, these features include (i) temporal attributes (day of week, day of month, month, year), (ii) clinician-specific statistics (historical availability rates, rolling weekly and monthly availability, group membership), and (iii) note-derived signals extracted from free-text scheduling comments (e.g., "Paid-time-off," Vacation," "Vacation," "Conference," "OR coverage"). Let $Y$ denote the observed outcomes, i.e., binary clinician-day availability labels indicating whether a clinician was available (1) or not (0). These labels are refined using unstructured notes to capture exceptions such as illness, paid-time-off, or cross-site redeployment. The predictor $h \in H$ is trained by minimizing a loss function $L$, mapping features $X$ to observed availability outcomes $Y$:

$$h^* = \arg\min_{h \in H} L(h(X), Y). \quad (1)$$

The resulting predictor $h^*$ produces estimates of uncertain parameters, which are then passed to the optimization step. Let $z$ denote the decision variables, and let $f(\cdot)$ denote the objective function. The optimization problem is formulated as:

$$z^* = \arg\max_{z} f(z, h^*(X), X), \quad (2)$$

subject to operational constraints. In this clinician scheduling application, $h^*(X)$ represents predicted clinician-day availability (probability scores), while $z^*$ represents the optimal assignments of clinicians to shifts.

Although PTO carries the risk that prediction errors may degrade optimization outcomes [23-25], problem-specific formulations can mitigate this sensitivity. In the proposed PTO-CS framework, this risk is mitigated by predicting probabilistic clinician–day availability scores rather than exact shift assignments. Probabilistic inputs are more robust to uncertainty and directly reflect clinicians' key priorities: compliance with contractual cFTE requirements, equitable workload distribution, and fair access to diverse shift types.

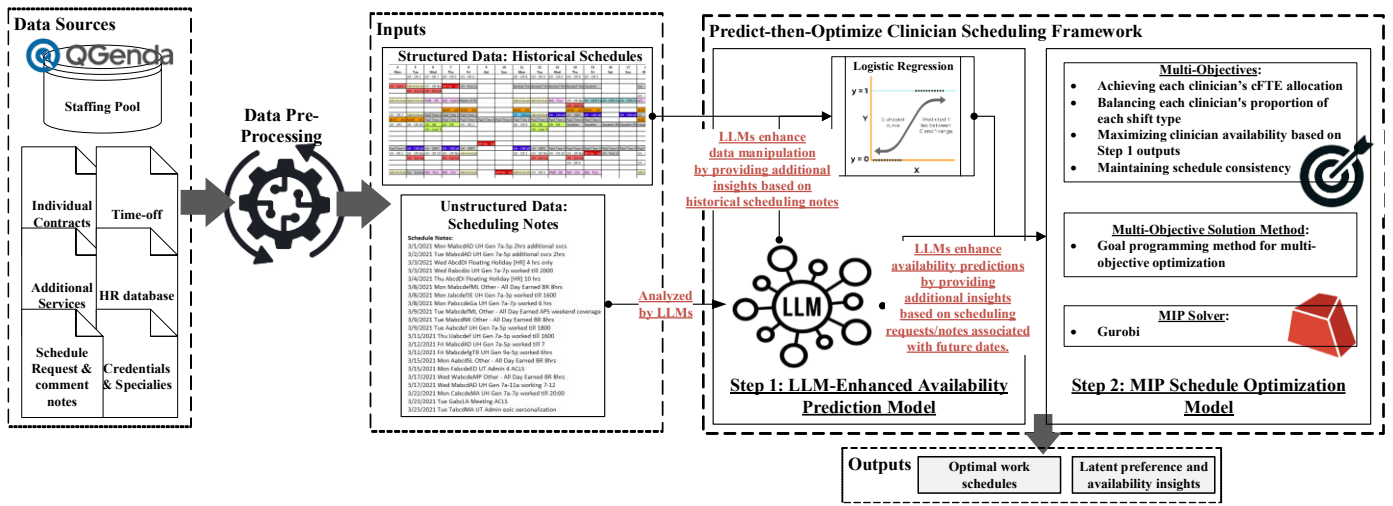

Fig. 2. The predict-then-optimize clinician scheduling (PTO-CS) framework

A major advantage of PTO is its modularity. Prediction and optimization modules can be validated independently, allowing stakeholders to verify intermediate outputs. This modular separation also ensures adaptability, i.e., advances in ML, particularly LLMs, can refine preference extraction, while new developments in optimization can improve the performance of the solvers and/or fairness formulations. In general, these features provide the flexibility, transparency, and interpretability required for adoption in real-world clinical scheduling practices.

### III. Research Methodology

The proposed PTO-CS framework (see Fig. 2) follows a PTO paradigm with two steps. In Step 1 (prediction), structured scheduling data and unstructured notes are transformed into calibrated clinician-day availability probabilities. In Step 2 (optimization), these probabilities are embedded into a multi-objective MIP that (i) ensures compliance with cFTE allocations, (ii) enforces equitable proportions of shift types, (iii) maximizes clinician availability for assigned shifts, and (iv) maintains schedule consistency. This modular design separates learning from decision-making, ensuring transparency and adaptability.

#### A. Data Simulation and Preprocessing

The dataset for this study was constructed from synthetic schedule templates enriched with LLM-generated notes, covering the period March 2021 through September 2024, to emulate the scheduling environment of the outpatient pain clinic.

Each month, an anonymous schedule template was produced by a month-by-month optimization model developed under requirements specified by the clinical informatics team. Clinician identifiers were masked as clinician IDs to ensure de-identification. The templates encoded institutional, operational, and individualized requirements supplied by the scheduling team, including coverage rules, workload distributions, and clinician-specific preferences. Designed to be compatible with commonly used scheduling platforms (e.g., Qgenda), these templates were also embedded in the scheduling processes at the concerned outpatient clinic. As a result, the templates provide a rich representation of real-world scheduling practices, balancing institutional/clinical policies with clinician-level preferences and constraints.

Although actual operations often introduce modifications (e.g., last-minute changes, urgent swaps, part-time call-ins), in this study, the monthly templates serve as proxies for historical schedules, providing a consistent clinician-day grid that reflects institutional intent. Each scheduled clinician-day pair is labeled as available, while non-scheduled pairs are labeled as not available after excluding weekends, holidays, and site-specific blackout days.

To simulate the unstructured component of scheduling data, LLMs external to the PTO-CS framework generated synthetic free-text comments. 10% of assigned shifts were randomly selected and annotated with plausible, contextually relevant notes (e.g., "paid-time-off," "Vacation," "Interview Day," "Conference," "Covering OR"), along with more natural, free-form requests and remarks consistent with clinical practice. This sampling rate reflects the typical frequency of annotation in operational scheduling systems, ensuring variability while preserving realism. By design, these notes inject implicit constraints and contextual information that structured templates alone do not capture.

From the structured schedules, we engineered temporal and clinician-specific features, including day of week, day of month, month, year, weekly availability rate, monthly availability rate, and weekday-specific availability rates. These features were merged with availability labels and simulated note annotations to form a consolidated modeling dataset, with one row per clinician-day.

The resulting dataset integrates synthetic scheduling data consisting of structured schedules and unstructured free-text scheduling notes and comments, establishing a realistic and auditable foundation for predictive modeling and downstream optimization. This dataset ensures alignment with operational practice and supports the practical applicability of the proposed PTO-CS framework.

*B. Step 1: LLM-Enhanced Availability Prediction Model*

Step 1 of the PTO-CS framework is a classification-based clinician availability prediction model enhanced with LLMs. While structured historical schedules capture the shift assignments, they do not always reflect actual availability. One-day exceptions, e.g., paid-time-off, illness, interview days, or cross-site coverage, appear only in unstructured free-text scheduling notes. Ignoring these notes introduces label noise, leading to overstated availability on conflict days and understated availability on others. Such miscalibration propagates into the prediction stage, undermining fairness, coverage, and workload balance in downstream scheduling.

To address this limitation, we combine a supervised, classification-based availability prediction model trained on structured features with an LLM-based note reader. The LLM strengthens availability prediction in two complementary ways: (i) by refining historical training labels through the extraction of conflicts and preferences from past notes, and (ii) by correcting future predictions when advance scheduling requests explicitly signal availability or unavailability. Together, these components ensure that predictions incorporate both structured scheduling patterns and contextual constraints, providing cleaner, more reliable inputs to the optimization stage.

*1) LLM-Driven Labels*

Unstructured scheduling notes often contain the real reason a clinician can or cannot work on a given day (e.g., paid-time-off, illness, conference, cross-site coverage, OR duty). To capture this information systematically, we employ a small, locally deployed instruction-tuned language model (Google FLAN-T5 base [26]) to convert notes into binary, auditable availability signals. These signals serve two purposes: (i) improving training labels for historical data, and (ii) refining probability estimates for upcoming schedules. FLAN-T5 is chosen because it handles short instruction-style prompts reliably, runs efficiently on standard hardware, and produces stable outputs without reliance on external APIs. Running locally also ensures privacy, reproducibility, and low cost. Given that notes are brief and map naturally to a two-label space, larger models would add little benefit relative to their higher computational costs. FLAN-T5 offers a balanced trade-off between accuracy, latency, and interpretability.

Formally, let $i$ index clinician and $t$ index the date. From the structured schedule, we extract a binary label $y_{i,t}^{struct} \in \{0,1\}$, where 1=available and 0=not available. When a schedule note exists for $(i,t)$, the LLM outputs a binary signal $l_{i,t} \in \{0,1\}$ with the same semantics. Overtime-related comments (e.g., "Overtime," "stayed late," "extended clinic") are mapped to $l_{i,t} = 1$, while conflicts (e.g., paid-time-off, illness, conference, or cross-site duties) are mapped to $l_{i,t} = 0$. If multiple comments exist, explicit conflicts take precedence; otherwise, the most recent non-conflict note is retained. All inferences run locally, and each decision is logged in a machine-readable audit trail that records the normalized text, matched rationale, and final label.

The final alignment signal is obtained by combining schedule and comment labels with a conservative rule as follows:

$$y_{i,t}^{\text{final}} = \begin{cases} y_{i,t}^{struct} \cdot l_{i,t}, & \text{if a comment exists for } (i,t), \\ y_{i,t}^{struct}, & \text{otherwise.} \end{cases} \quad (3)$$

This ensures that explicit unavailability overrides the schedule, while neutral or positive notes leave the original schedule unchanged. In the subsequent probability model, we apply the same conservative rule, i.e., probabilities are hard zeroed when comments indicate unavailability; otherwise, the calibrated baseline is preserved. This rule reduces false positives on conflict days, maintains calibrated nuance elsewhere, and integrates seamlessly into existing scheduling workflows without altering clinic rules.

*2) Classification-based Models*

Before integrating LLM-derived labels, we benchmarked supervised classifiers trained solely on structured scheduling features. These features, engineered from common scheduling platform records, include temporal attributes (day-of-week, day-of-month, month, year) and clinician-specific statistics (historical availability, rolling weekly and monthly availability rates, and clinician-group indicators). The goal is to estimate the probability that a clinician $i$ would be available on a given day $t$, producing calibrated probabilities suitable for downstream work schedule optimization.

Our primary estimator is logistic regression, selected for its interpretability and probabilistic outputs. Formally, let $a_{i,t} \in \{0,1\}$ denote the day-level availability indicator and $X_{i,t} \in \mathbb{R}^K$ the engineered feature vector with $K$ features. The probability of availability from the logistic regression classifier is formulated as follows:

$$p_{i,t}^C = Pr(a_{i,t} = 1 | X_{i,t}) = \sigma(\beta_0 + \beta^T X_{i,t}) = \frac{1}{1+\exp\{-(\beta_0 + \beta^T X_{i,t})\}}, \quad (4)$$

where $\sigma(\eta) = \frac{1}{1+e^{-\eta}}$ is the logistic function, $\beta_0$ is the intercept, and $\beta \in \mathbb{R}^K$ is the coefficient vector. This benchmark evaluation respected temporal order: using forward-chaining validation [27], we trained on earlier months and validated on later ones, following how scheduling decisions are made in practice. Class imbalance across the full clinician-day grid was addressed through resampling, and probability calibration ensured that predicted probabilities corresponded to actual availability rates. Specifically, month-ahead calibration was performed by learning a mapping between raw predictions and observed outcomes in the immediately following month. For each month, predicted probabilities were grouped into bins, the empirical availability frequency within each bin was computed, and these frequencies were then used to adjust subsequent predictions. This rolling calibration procedure avoided temporal leakage while aligning predicted probabilities with real-world availability patterns. Logistic regression achieved strong operational performance, with Accuracy=0.78, Recall for the available class=0.85, and Macro-F1=0.78 (Table I). The high recall reduced missed staffing opportunities, and the calibrated probabilities were stable month-to-month, making this model both reliable and auditable for clinical partners.

To benchmark against nonlinear alternatives, we also evaluated decision trees and random forests. The decision tree model recursively partitioned the feature space to capture threshold and interaction effects that a linear model might miss (e.g., "Fridays in July for Group A"). With controlled depth and pruning by default, the tree achieved Accuracy=0.79, Precision for the available class=0.83, Recall=0.75, and Macro-F=0.79. Although interpretable, single trees were less stable across months and required post-hoc calibration for probability estimates.

The random forest model, which averages many decorrelated trees, achieved the highest overall Accuracy=0.80 and Macro-F1=0.79. It generalized well across clinicians and seasons, but its raw probability outputs were miscalibrated, and the model is not as interpretable as logistic regression or single trees.

In summary, structured-only baselines reflected the classic trade-off between interpretability and accuracy. Logistic regression offered competitive accuracy with the strongest recall for the available class, decision trees captured nonlinear rules at the cost of stability, and random forests improved overall accuracy while reducing transparency. Taken together, these results supported the use of logistic regression as the default baseline $p_{i,t}^C$.

In the next subsection, we extend this baseline by integrating note-derived signals from the LLM, ensuring that explicit conflicts override structured predictions and yielding a single auditable probability for Step 2 of the PTO-CS framework.

*3) Integrating LLM Models and Classifiers*

Step 1 of the PTO-CS framework culminates in an LLM–enhanced availability prediction model that integrates structured classification-based probabilities with contextual corrections derived from unstructured scheduling notes. The process unfolded in three sequential stages:

*Temporal Segregation*: LLM-extracted labels from past notes are merged with structured schedule features (e.g., day-of-week, rolling availability, clinician group) to generate refined training data. This integration reduces label noise introduced when exceptions such as paid-time-off or cross-site coverage are ignored, yielding more accurate historical availability baselines.

*Probabilistic Estimation*: Classification models, led by logistic regression with calibrated outputs, generate baseline probabilities of clinician availability for future days. E.g., the model may estimate that clinician $i$ has a 92% probability of being available on Week 2, Day 3, based on temporal and historical scheduling patterns, which serves as the baseline availability probability.

*LLM-Guided Correction*: For future schedules, clinicians often submit notes alongside their requests (e.g., "Interview Day," "Vacation," "Conference," or "OR coverage"). The LLM processes these notes and, when they explicitly indicate unavailability, the corresponding probability is hard-set to zero. Otherwise, the baseline probability is retained. This conservative adjustment enforces deterministic corrections where textual evidence is clear while preserving probabilistic nuance elsewhere.

This three-stage design ensures that the model retains probabilistic richness when structured features and notes align, while enforcing strict corrections when conflicts are explicitly documented. By combining the structured regularities captured by logistic regression with the contextual awareness of FLAN-T5, the LLM–enhanced availability prediction model balances interpretability, adaptability, and predictive accuracy, resulting in availability estimates that are robust, transparent, and operationally actionable.

Formally, the final probability is refined conservatively as follows:

$$p_{i,t} = \begin{cases} 0, & \text{if a comment exists and } l_{i,t} = 0, \\ p_{i,t}^C, & \text{otherwise.} \end{cases} \quad (5)$$

This formulation hard-zeros availability when notes clearly indicate a conflict, while leaving the calibrated baseline unchanged in all other cases. The result is a single, auditable probability estimate that reduces false positives on conflict days, preserves calibrated probabilities elsewhere, and integrates seamlessly into existing scheduling workflows without altering institutional rules.

## C. Step 2: Optimization-Based Schedule Generation

Step 2 of the PTO-CS framework embeds the LLM–enhanced availability predictions from Step 1 into a multi-objective MIP model that produces operational schedules. While Step 1 refines probabilistic estimates of clinician availability, Step 2 transforms these signals into optimized schedules that account for institutional requirements, workload balance, and clinician-level fairness.

*1) Decision Variables*

Let $I$ denote the set of clinicians, $S$ the set of shift types, and $D$ the set of days in the scheduling horizon. The binary decision variable is defined as:

$$B_{i,s,t} = \begin{cases} 1, & \text{if } i \in I \text{ is assigned to } s \in S \text{ on } t \in D, \\ 0, & \text{otherwise.} \end{cases} \quad (6)$$

Each clinician-day pair is associated with a probabilistic availability estimate $p_{i,t} \in [0,1]$ derived in Step 1. This probability is not a decision variable but informs assignment feasibility and objective weighting.

*2) Objectives*

Step 2 model considers four objectives that are formulated as follows:

$$z_1 = \max_B \{-\sum_{i \in I} |\sum_{s \in S} \sum_{t \in D} B_{i,s,t} - cFTE_i|\}, \quad (7)$$

$$z_2 = \max_B \{-\sum_{i \in I} \sum_{s \in S} \left|\frac{1}{|D|} \sum_{t \in D} B_{i,s,t} - \overline{B}_s\right|\}, \quad (8)$$

$$z_3 = \max_B \{-\sum_{i \in I} \sum_{s \in S} \sum_{t \in D} p_{i,t} B_{i,s,t}\}, \quad (9)$$

TABLE I. ML MODEL PERFORMANCE BENCHMARK

| Model | Accuracy | Precision (0) | Precision (1) | Recall (0) | Recall (1) | F1 (0) | F1 (1) | Macro F1 |
|---|---|---|---|---|---|---|---|---|
| Logistic Regression | 0.78 | 0.82 | 0.76 | 0.71 | **0.85** | 0.76 | 0.80 | 0.78 |
| Random Forest | **0.80** | 0.79 | 0.80 | 0.79 | 0.80 | 0.79 | 0.80 | **0.79** |
| Decision Tree | 0.79 | 0.76 | **0.83** | **0.84** | 0.75 | **0.80** | 0.79 | 0.79 |

$$z_4 = \max_B \{-\sum_{i\in I}\sum_{s\in S}\sum_{t\in D}|B_{i,s,t} - B_{i,s,t}^{prev}|\}. \quad (10)$$

Objective (7) enforces compliance with contractual cFTEs by penalizing deviations between the assigned workload and the target allocation. Objective (8) promotes equity by balancing the distribution of shift types across clinicians. For each shift type $s \in S$, the deviation of a clinician's assignment frequency from the average share $\bar{B}_s$ is minimized. Objective (9) maximizes alignment between assignments and predicted availability. Objective (10) seeks schedule consistency by minimizing disruptive changes relative to the prior schedule $B_{i,s,t}^{prev}$.

*3) Constraints*

Constraints of Step 2 model are listed in (11)–(15) in the following:

$$\sum_{i\in I} B_{i,s,t} = r_{s,t}, \forall s \in S, t \in D, \quad (11)$$

$$\sum_{s\in S} B_{i,s,t} \leq 1, \forall i \in I, t \in D, \quad (12)$$

$$Li \leq \sum_{s\in S}\sum_{t\in D} B_{i,s,t} \leq Ui, \forall i \in I \quad (13)$$

$$B_{i,s,t} \leq \bar{a}_{i,t}, \forall i \in I, s \in S, t \in D \quad (14)$$

$$\sum_{s\in S}\sum_{t\in D^{wkend}} B_{i,s,t} \leq W_i, \forall i \in I \quad (15)$$

Constraints (11) ensure that the total number of clinicians assigned to each shift matches demand, where $r_{s,t}$ is the required staffing level. Constraints (12) ensure each clinician may work at most one shift per day. Constraints (13) are the workload bounds constraints, which ensure the total number of shifts assigned to each clinician is within a predefined range $[Li, Ui]$, where $Li$ and $Ui$ are the lower and upper bounds for clinician $i$. These bounds enforce fairness and ensure compliance with contractual workload agreements. Constraints (14) ensure that assignments must also respect hard availability restrictions derived from Step 1, where $\bar{a}_{i,t} \in \{0,1\}$ is the hard availability indicator. Constraints (15) ensure institution-specific policies, e.g., limits on weekend or night duties, are encoded, where $D^{wkend}$ is the set of weekend days and $W_i$ is the allowable weekend workload for clinician $i$.

*4) Solution Method: Lexicographic Goal Programming*

The optimization problem in Step 2 integrates four conflicting objectives. Because these objectives cannot generally be optimized simultaneously without trade-offs, a lexicographic goal programming approach is adopted to balance competing requirements.

In lexicographic goal programming, each objective is reformulated as a goal with a specified aspiration level $g_k$. Deviational variables $d_k^-$ and $d_k^+$ measure under- and over-achievement relative to the target. The model then minimizes a weighted sum of these deviations:

$$Z = \min \sum_{k=1}^{4} \lambda_k(d_k^- + d_k^+), \quad (16)$$

s.t. $z_k + d_k^- - d_k^+ = g_k, d_k^-, d_k^+ \geq 0, \forall k \in \{1,\dots,4\}, \quad (17)$

where $\lambda_k$ represents the relative priority weight assigned to each goal. In the PTO-CS framework, the priorities are defined as:

*Goal 1*: cFTE Compliance. Ensuring that each clinician's total number of assigned shifts matches their contractual cFTE level is treated as the highest-priority requirement. Any deviation from this allocation is minimized before addressing other objectives.

*Goal 2*: Equitable Shift Type Proportions. Given cFTE compliance, the model seeks to distribute different shift types equitably across clinicians. Deviations capture imbalances, particularly in undesirable or high-burden shifts.

*Goal 3*: Availability Maximization. Subject to Goals 1 and 2, schedules are optimized to maximize alignment with predicted availability probabilities $p_{i,t}$, reducing the likelihood of infeasible or last-minute changes.

*Goal 4*: Schedule Consistency. Finally, consistency with the previous month's schedule is enforced to the extent possible, minimizing disruptive changes while respecting higher-priority goals.

This lexicographic structure reflects clinical and institutional priorities: contractual cFTE obligations must be satisfied first, equitable distribution of shift types is preserved to reduce workload disparities, and predicted availability is leveraged to mitigate uncertainty in clinical resources. Consistency with prior schedules is incorporated as a secondary objective, recognizing its value for clinician satisfaction while ensuring it does not override fairness or compliance. By adopting this formulation, PTO-CS produces schedules that are feasible, policy-compliant, equitable, and robust to variability in clinician availability, while limiting unnecessary disruptions to clinicians' work patterns.

IV. EVALUATION AND RESULTS

*A. Evaluation Setup and Metrics*

We evaluated the proposed PTO-CS framework using realistic synthetic outpatient pain clinic scheduling data with schedule comments as described in Section III. For evaluation purposes, this dataset was treated as historical (actual) schedules, against which optimized schedules generated by PTO-CS were compared.

The evaluation period spanned six consecutive months, i.e., from March 2024 to August 2024. This timeframe was selected for two main reasons. First, it contains seasonal diversity in scheduling patterns, including the transition from the academic spring semester into summer. During this period, clinician availability is affected by increased clinical duties, vacations, and professional conferences. Importantly, as part of an academic anesthesiology department, schedules are also influenced by new residents beginning rotations in July. More experienced faculty, often holding dual appointments in both the OR and the outpatient clinic, are more likely to be redeployed to high-acuity sites (e.g., ORs) when resident staffing is relatively less experienced. Second, this window offered stability in clinical scheduling operations, as scheduling processes and templates had been standardized earlier in the year, ensuring consistent data quality since the synthetic data is based on historical schedule templates. Overall, this six-month window reflects a mix of routine and high-variability conditions while operations remained stable, providing a natural stress test for PTO-CS. Section IV-C and Fig. 3 report the month-by-month results.

The evaluation followed a two-level design. At the distributional level, we compared clinician counts per shift type on each and across the schedule horizon for both optimized and historical schedules. The analysis emphasized two core shift types, i.e., Clinic (Clin) and Procedure (Proc), which served as proxies for clinical coverage compliance and cFTE allocation fairness. At the instance level, we conducted day-by-day clinician–day comparisons. An assignment was considered a match if the same clinician was scheduled to work on that day in both the optimized and historical schedules, regardless of the specific shift type; otherwise, it was treated as a mismatch. Match accuracy, expressed as the percentage of clinician–day matches, was used as an indicator of operational alignment rather than predictive accuracy in the machine learning sense.

Optimization experiments were conducted in Python 3.13, with models solved using Gurobi 11.0 on a workstation (Intel Core i7-14700, 32 GB RAM). Prediction baselines (logistic regression, decision tree, random forest) and the LLM-enhanced availability predictor were implemented in Python (scikit-learn). The note-classification module employed a locally deployed Google FLAN-T5 base model, ensuring reproducibility, privacy-preserving inference, and seamless prediction–optimization integration.

### B. Key Performance Indicators (KPIs)

*Coverage rate.* The clinical coverage rate of required shifts that were filled in the schedule:

$$Coverage\ rate = \frac{filled\ shifts}{Total\ required\ shifts} \times 100\%. \quad (18)$$

*Total cFTE deviation.* The total absolute difference of each clinician's assigned workload, i.e., shifts assigned in the outpatient clinic divided by the number of duty days in the month, and their contractual cFTE. For clinician $i$:

$$Total\ cFTE\ deviation = \sum_{i \in I} cFTE\ deviation_i, \quad (19)$$

where,

$$cFTE\ deviation_i = \left| \frac{Assigned\_clinic\_shifts_i}{Duties\ days\ in\ month} - cFTE_i \right|. \quad (20)$$

*Shift-type proportional equity index.* The fairness indicator of distributing shift type $s$ across clinicians. For each shift type $s$, we computed the clinician-level proportions:

$$Equity\ index_i = \sum_{s \in S} Var(c_{i,s}), \quad (21)$$

where $c_{i,s}$ denotes the proportion of clinician $i$ assigned to shift $s$, formally:

$$c_{i,s} = \frac{Shifts\ of\ type\ s\ assigned\ to\ clinician\ i}{Total\ shifts\ assigned\ to\ clinician\ i}. \quad (22)$$

*Match accuracy.* The proportion of clinician–day assignments in the optimized schedule that matched the historical schedule:

$$Match\ Accuracy = \frac{Matched\ assignments}{Total\ assignments} \times 100\%. \quad (23)$$

### C. Results and Interpretation

The performance of PTO-CS is evaluated using the four KPIs described above and compared against historical scheduling practices over a six-month horizon (March–August 2024). Monthly results were consistent across the window.

Coverage rate performance is summarized in Fig. 3. Historical schedules showed incomplete coverage in multiple months, with Clin shift coverage ranging from 0.77 (April) to 0.97 (March) and Proc shift coverage ranging from 0.68 (August) to 1.0 (April). In contrast, PTO-CS consistently achieved 100% coverage for both Clin and Proc shifts across all six months. These results confirm that the optimization framework effectively enforced coverage constraints (Eq. (11)), ensuring all required shifts were filled regardless of seasonal fluctuations in clinician availability.

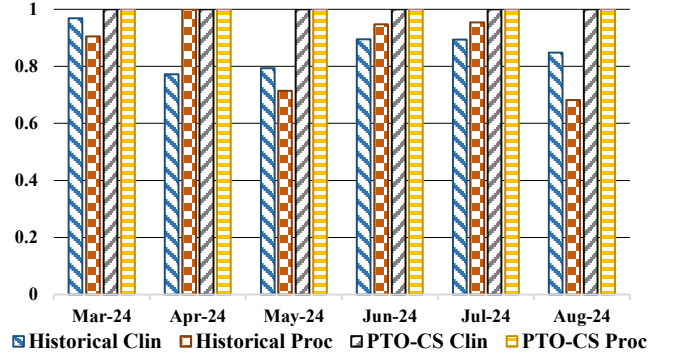

Fig. 3. Coverage rate, historical schedule versus PTO-CS schedule.

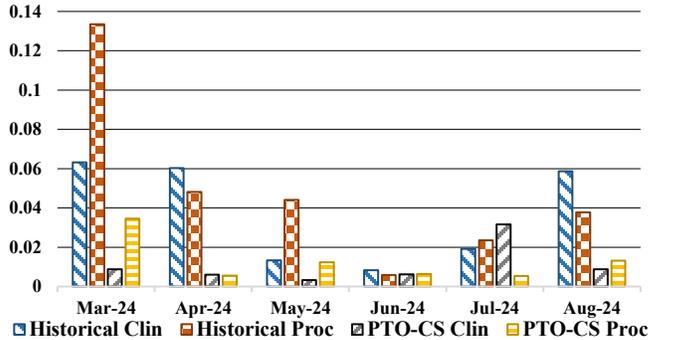

Fig. 4. Variances of proportions of shift types, historical schedules versus PTO-CS schedules.

Shift-type proportional equity results are summarized in Fig. 4. Historical schedules exhibited higher variance in clinician-level allocations, indicating workload imbalances across both Clin and Proc shifts. E.g., variance in Proc assignments reached 0.134 in March and remained elevated in April (0.048) and May (0.044). By contrast, PTO-CS consistently reduced variance, with Proc shift proportions ranging from 0.005 to 0.034 across all months. Clinic shift allocations showed a similar pattern, i.e., historical variance peaked at 0.063 in March, whereas PTO-CS reduced this to below 0.009 in most months.

Overall, PTO-CS produced substantially more balanced distributions of both Clin and Proc shifts across clinicians. These improvements demonstrate that the framework effectively addressed proportional equity objectives, ensuring fairer allocation of workload compared with historical scheduling practices.

cFTE compliance and schedule accuracy results are shown in Fig. 5. Historical schedules demonstrated substantial misalignment with contractual cFTEs, with total cFTE deviations ranging from 0.45 (March) to 1.25 (April). In contrast, PTO-CS reduced misalignment to between 0.15 and 0.32, representing a marked improvement in aligning clinician assignments with contractual obligations.

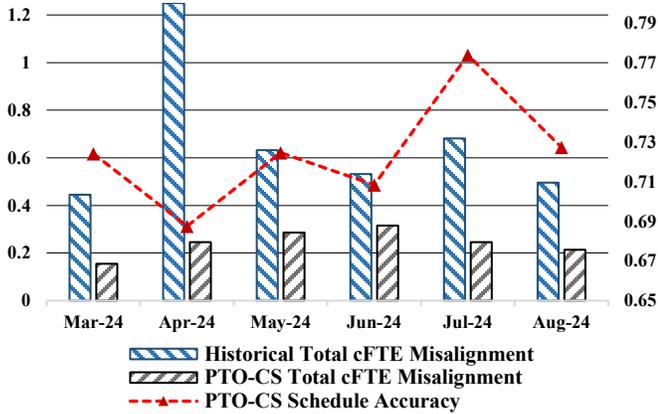

Fig. 5. PTO-CS accuracy and cFTE misalignment comparison

Schedule accuracy for PTO-CS, measured as clinician–day alignment with historical schedules, ranged from 68.8% (April) to 77.4% (July). Accuracy was lower in months where fairness and coverage constraints were binding (e.g., April and June) but improved when clinician workloads were more flexible. These results illustrate the intended trade-off in the lexicographic optimization, i.e., reductions in match accuracy were offset by significant improvements in fairness and compliance with contractual cFTE levels.

## V. CONCLUSION AND FUTURE WORK

In this work, PTO-CS, a predict-then-optimize framework for outpatient clinician scheduling, is developed and evaluated. By combining LLM-enhanced availability predictions with multi-objective optimization, the PTO-CS achieves full coverage, improves workload fairness by enforcing an equitable proportion of shift types, and ensures compliance with contractual cFTE allocations benchmarked against historical schedules. The lower assignment match accuracy, i.e., comparing PTO-CS generated schedules with realistic historical schedule templates, reflects intentional reallocations to prioritize fairness and compliance, highlighting the trade-offs inherent in lexicographic optimization. The PTO-CS framework builds on prior consultations with the concerned outpatient pain clinic and its scheduling team, ensuring alignment with day-to-day operations and established workflows. This grounding in clinical practice enhances ease of adoption and underscores its potential for real-world deployment and measurable impact on clinical operations.

Several important limitations highlight directions for future work. First, the current application of LLMs was intentionally conservative, limited to extracting explicit conflicts (e.g., paid-time-off, illness, cross-site coverage) and mapping them to availability restrictions. While this approach improved baseline predictions, it underutilized the broader potential of LLMs. In addition, the fine-tuned FLAN-T5 model was selected for its efficiency, ease of deployment on standard hardware, and ability to run locally (addressing data safety concerns), which made it a practical choice for early-stage testing and proof of concept. Building on this foundation, future work will investigate alternative LLM architectures and fine-tuning strategies to improve the accuracy of extracting conflicts and overtime indicators. Furthermore, expanding the role of LLMs to capture nuanced preferences, partial availability (e.g., availability for two hours beyond a regular shift), and contextual signals beyond binary overrides represents a promising avenue for improving the practical value of this work.

Second, the PTO-CS is evaluated on realistic, close-to-final schedule templates augmented with synthetic LLM-generated notes. Future work will proceed along two complementary directions. On one hand, we will develop a rigorous LLM-based schedule note simulation tool to better mimic the density and variability of real scheduling notes. Since clinicians may or may not choose to leave notes, this tool will allow us to systematically assess how the performance of PTO-CS varies under different levels of note density. Additionally, we will establish benchmarks for LLM-based conflict extraction to evaluate reliability and quantify misclassification risks. On the other hand, access to real scheduling records will be essential to measure the framework's practical impact on clinical operations. We are actively collaborating with clinical partners to adapt PTO-CS to operational workflows and to validate it prospectively on authentic scheduling documents and notes.

Third, future work will also expand schedule performance evaluation beyond operational metrics, e.g., coverage, fairness, and workload compliance, to include downstream outcomes, e.g., patient access, quality of care, clinician well-being, and fatigue. While the six-month evaluation window captured both routine and high-variability conditions, longer-term studies are needed to assess cumulative impacts, e.g., fatigue and burnout.

Methodologically, the current optimization relied on a fixed lexicographic priority structure. Exploring alternative multi-objective solution approaches, e.g., ε-constraint methods or interactive procedures [28], could provide decision makers with greater flexibility. Moreover, the decoupled predict-then-optimize paradigm introduces risks of error propagation and misalignment between prediction and decision-making. Addressing these challenges through integrated joint prediction-and-optimization approaches [29] that integrate predictive modeling with scheduling objectives represents an important future direction.

In summary, PTO-CS offers a transparent, modular, and scalable approach to clinician scheduling that balances fairness, coverage, and workload compliance. With enhanced use of LLMs, real-world validation, and methodological refinements, the framework has the potential to reduce scheduler burden, improve workload equity, and support sustainable clinical operations in diverse healthcare settings.

ACKNOWLEDGMENT

The authors thank the faculty and staff of the Department of Anesthesiology at the University of Texas Health Science Center at San Antonio for their valuable feedback in the development of this work. The authors also acknowledge the use of OpenAI's ChatGPT to assist with language polishing of the manuscript.

REFERENCES

[1] L. H. Aiken, S. P. Clarke, D. M. Sloane, J. Sochalski, and J. H. Silber, "Hospital nurse staffing and patient mortality, nurse burnout, and job dissatisfaction," *JAMA*, vol. 288, no. 16, p. 1987, Oct. 2002, doi: 10.1001/jama.288.16.1987.

[2] Y. Lan, D. Goradia, and A. Chandrasekaran, "Ancillary cost implications of physicians multisiting and inter-organizational collaboration during healthcare delivery," *Prod. Oper. Manag.*, vol. 31, no. 2, pp. 561–582, 2022, doi: 10.1111/poms.13567.

[3] "AAMC Report Reinforces Mounting Physician Shortage," AAMC. Accessed: May 09, 2025. [Online]. Available: https://www.aamc.org/news/press-releases/aamc-report-reinforces-mounting-physician-shortage

[4] A. Al Yahyaei, A. Hewison, N. Efstathiou, and D. Carrick-Sen, "Nurses' intention to stay in the work environment in acute healthcare: a systematic review," *J. Res. Nurs.*, vol. 27, no. 4, pp. 374–397, June 2022, doi: 10.1177/17449871221080731.

[5] Ringo, "What healthcare CFOs should know about the CRNA shortage." Accessed: May 09, 2025. [Online]. Available: https://www.goringo.com/blog/what-healthcare-cfos-should-know-about-the-crna-shortage

[6] S. Rath, K. Rajaram, and A. Mahajan, "Integrated anesthesiologist and room scheduling for surgeries: Methodology and application," *Oper. Res.*, vol. 65, no. 6, p. 35, 2017.

[7] T. Cayirli, E. Veral, and H. Rosen, "Designing appointment scheduling systems for ambulatory care services," *Health Care Manag. Sci.*, vol. 9, no. 1, pp. 47–58, Feb. 2006, doi: 10.1007/s10729-006-6279-5.

[8] B. Mets, "A career in academic anesthesiology," in *A guide to Anesthesiology for medical students*, 2018. Accessed: May 09, 2025. [Online]. Available: https://www.asahq.org/-/media/sites/asahq/files/public/about-asa/component-societies/med-students/medical-student-guide/chapter_5_guide-to-a-career-in-anesthesiology-for-medical-students.pdf

[9] S. Hoffman and I. Ganguli, "Access to Primary Care and Health Care Fragmentation," *Fac. Publ.*, Jan. 2026, [Online]. Available: https://scholarlycommons.law.case.edu/faculty_publications/2327

[10] D. Gupta and B. Denton, "Appointment scheduling in health care: Challenges and opportunities," *IIE Trans.*, vol. 40, no. 9, pp. 800–819, July 2008, doi: 10.1080/07408170802165880.

[11] A. Erekat, G. Servis, S. C. Madathil, and M. T. Khasawneh, "Efficient operating room planning using an ensemble learning approach to predict surgery cancellations," *IISE Trans. Healthc. Syst. Eng.*, vol. 10, no. 1, pp. 18–32, Jan. 2020, doi: 10.1080/24725579.2019.1641576.

[12] S. Youn, H. N. Geismar, C. Sriskandarajah, and V. Tiwari, "Adaptive Capacity Planning for Ambulatory Surgery Centers," *Manuf. Serv. Oper. Manag.*, May 2022, doi: 10.1287/msom.2022.1109.

[13] A. Ahmadi-Javid, P. Seyedi, and S. S. Syam, "A survey of healthcare facility location," *Comput. Oper. Res.*, vol. 79, pp. 223–263, Mar. 2017, doi: 10.1016/j.cor.2016.05.018.

[14] A. Azari, V. P. Janeja, and S. Levin, "Imbalanced learning to predict long stay Emergency Department patients," in *2015 IEEE International Conference on Bioinformatics and Biomedicine (BIBM)*, Nov. 2015, pp. 807–814. doi: 10.1109/BIBM.2015.7359790.

[15] Y. Hu *et al.*, "Use of Real-Time Information to Predict Future Arrivals in the Emergency Department," *Ann. Emerg. Med.*, vol. 81, no. 6, pp. 728–737, June 2023, doi: 10.1016/j.annemergmed.2022.11.005.

[16] K. Sun, A. Roy, and J. M. Tobin, "Artificial intelligence and machine learning: Definition of terms and current concepts in critical care research," *J. Crit. Care*, vol. 82, p. 154792, Aug. 2024, doi: 10.1016/j.jcrc.2024.154792.

[17] R. J. Tobin, C. R. Walker, R. Moss, J. M. McCaw, D. J. Price, and F. M. Shearer, "A modular approach to forecasting COVID-19 hospital bed occupancy," *Commun. Med.*, vol. 5, no. 1, p. 349, Aug. 2025, doi: 10.1038/s43856-025-01086-0.

[18] A. Ramamurthi *et al.*, "Applying Large Language Models for Surgical Case Length Prediction," *JAMA Surg.*, vol. 160, no. 8, pp. 894–902, Aug. 2025, doi: 10.1001/jamasurg.2025.2154.

[19] L. Tang *et al.*, "Evaluating large language models on medical evidence summarization," *Npj Digit. Med.*, vol. 6, no. 1, p. 158, Aug. 2023, doi: 10.1038/s41746-023-00896-7.

[20] D. Van Veen *et al.*, "Adapted large language models can outperform medical experts in clinical text summarization," *Nat. Med.*, vol. 30, no. 4, pp. 1134–1142, Apr. 2024, doi: 10.1038/s41591-024-02855-5.

[21] A. Anis Lahoud, A. S. Khan, E. Schaffernicht, M. Trincavelli, and J. A. Stork, "Predict-and-Optimize Techniques for Data-Driven Optimization Problems: A Review," *Neural Process. Lett.*, vol. 57, no. 2, p. 40, Apr. 2025, doi: 10.1007/s11063-025-11746-w.

[22] K. Sun, M. Sun, D. Agrawal, R. Dravenstott, F. Rosinia, and A. Roy, "Equitable anesthesiologist scheduling under demand uncertainty using multiobjective programming," *Prod. Oper. Manag.*, vol. 32, no. 11, pp. 3699–3716, 2023, doi: 10.1111/poms.14058.

[23] P. Donti, B. Amos, and J. Z. Kolter, "Task-based End-to-end Model Learning in Stochastic Optimization," in *Advances in Neural Information Processing Systems*, Curran Associates, Inc., 2017. Accessed: Aug. 23, 2025. [Online]. Available: https://proceedings.neurips.cc/paper/2017/hash/3fc2c60b5782f641f76bcefc39fb2392-Abstract.html

[24] D. Grimes, G. Ifrim, B. O'Sullivan, and H. Simonis, "Analyzing the impact of electricity price forecasting on energy cost-aware scheduling," *Sustain. Comput. Inform. Syst.*, vol. 4, no. 4, pp. 276–291, Dec. 2014, doi: 10.1016/j.suscom.2014.08.009.

[25] G. Ifrim, B. O'Sullivan, and H. Simonis, "Properties of Energy-Price Forecasts for Scheduling," in *Principles and Practice of Constraint Programming*, M. Milano, Ed., Berlin, Heidelberg: Springer, 2012, pp. 957–972. doi: 10.1007/978-3-642-33558-7_68.

[26] J. Oza and H. Yadav, "Enhancing Question Prediction with Flan T5 -A Context-Aware Language Model Approach," 2023, *Authorea Preprints*. Accessed: Aug. 29, 2025. [Online]. Available: https://www.authorea.com/doi/full/10.22541/au.170258918.81486619?commit=9ee249a0d74897c207fdfeeee395f0819798f10f

[27] C. Bergmeir and J. M. Benítez, "On the use of cross-validation for time series predictor evaluation," *Inf. Sci.*, vol. 191, pp. 192–213, May 2012, doi: 10.1016/j.ins.2011.12.028.

[28] M. Sun, "Encyclopedia of Business Analytics and Optimization," in *Multiobjective Programming*, United States of America: Business Science Reference, 2014, pp. 1585–1604. doi: 10.4018/978-1-4666-5202-6.ch143.

[29] J. Dias Garcia, A. Street, T. Homem-de-Mello, and F. D. Muñoz, "Application-Driven Learning: A Closed-Loop Prediction and Optimization Approach Applied to Dynamic Reserves and Demand Forecasting," *Oper. Res.*, vol. 73, no. 1, pp. 22–39, Jan. 2025, doi: 10.1287/opre.2023.0565.